\documentclass[11pt]{article}
\usepackage[all]{xy}
\usepackage{amstex}
\usepackage{delarray}
\usepackage{amssymb}
\newtheorem{defi}{ Definition :}[section]

\newtheorem{theo}[defi]{Theorem }
\newtheorem{lem}[defi]{Lemma }
\newtheorem{prop}[defi]{Proposition }
\newtheorem{cor}[defi]{Corollary }
\newtheorem{rem}[defi]{Remark }
\newtheorem{rems}[defi]{Remarks }

\newcommand{\N}{{\Bbb N}}
\newcommand{\Z}{{\Bbb Z}}
\newcommand{\K}{{\Bbb K}}
\newcommand{\R}{{\Bbb R}}
\newcommand{\C}{{\Bbb C}}

\newcommand{\qh}{{quasihomogeneous }}
\newcommand{\degp}{{{\rm {N}}}}
\newcommand{\fct}{{{\cal F}(\K^n)}}
\newcommand{\chp}{{{\cal X}(\K^n)}}
\def\lcf{\lbrack\! \lbrack}
\def\rcf{\rbrack\! \rbrack}

\begin{document}
\footnote{{\bf {Key-words :}} Nambu-Poisson structures, singularities, Nambu-Poisson cohomology}
\footnote{{\bf {AMS classification :}} 53D17 }
\addtocounter{footnote}{-1}
\begin{center}
{\bf {\Large Computations of Nambu-Poisson cohomologies}}
\end{center}
\vspace{0.2 cm}
\begin{center}
Philippe MONNIER
\end{center}
\vspace{0.5 cm}
\begin{abstract}
In this paper, we want to associate to a $n$-vector on a manifold of dimension $n$ a cohomology
which generalizes the Poisson cohomology of a 2-dimensional Poisson manifold. Two possibilities are 
given here. One of them, the Nambu-Poisson cohomology, seems to be the most pertinent. We study 
these two cohomologies locally, in the case of germs of $n$-vectors on $\K^n$ 
($\K=\R {\mbox{ or }} \C$).
\end{abstract}

\section{Introduction}
A way to study a geometrical object is to associate to it a cohomology. In this paper, we focus
on the $n$-vectors on a $n$-dimensional manifold $M$. \\
If $n=2$, the 2-vectors on $M$ are the Poisson stuctures thus, we can consider the Poisson 
cohomology. In dimension 2, this cohomology has three spaces. The first one, $H^0$, is the 
space of functions whose Hamiltonian vector field is zero (Casimir functions).
The second one, $H^1$, is the quotient of the space of infinitesimal automorphisms 
(or Poisson vector fields) by the subspace of Hamiltonian vector fields. 
The last one, $H^2$, describes the deformations of the Poisson structure. 
In a previous paper ([Mo]), we have computed the cohomology of germs at 0 of
Poisson structures on $\K^2$ ($\K=\R {\mbox{ or }} \C$).\\
In order to generalize this cohomology to the $n$-dimensional case ($n\geq 3$), we can follow the
same reasoning. These spaces are not necessarily of finite dimension and it is not always
easy to describe them precisely.\\
\indent Recently, a team of Spanish researchers has defined a cohomology, called Nambu-Poisson
cohomology, for the Nambu-Poisson structures (see [I2]). In this paper, we adapt their 
construction
to our particular case. We will see that this cohomology generalizes in a certain sense the
Poisson cohomology in dimension 2. Then we compute locally this cohomology for germs at 0 of 
$n$-vectors
$\Lambda=f\frac{\partial }{\partial x_1}\wedge\hdots\wedge\frac{\partial }{\partial x_n}$ on
$\K^n$ ($\K=\R {\mbox { or }} \C$), with the assumption that $f$ is a \qh polynomial of finite
codimension (''most of'' the germs of $n$-vectors have this form). 
This computation is based on a preliminary result that we have shown, in the
formal case and in the analytical case (so, the ${\cal C}^\infty$ case is not entirely solved).
The techniques we use in 
this paper are quite the same as in [Mo].

\section{Nambu-Poisson cohomology}
Let $M$ be a differentiable manifold of dimension $n$ ($n\geq 3$), admitting a volume form
$\omega$. We denote ${\cal C}^\infty(M)$ the space 
of ${\cal C}^\infty$ functions on $M$, $\Omega^k(M)$ ($k=0,\hdots,n$) the 
${\cal C}^\infty(M)$-module of $k$-forms on $M$, and ${\cal X}^k(M)$ ($k=0,\hdots,n$) the  
${\cal C}^\infty(M)$-module of $k$-vectors on $M$.\\
We consider a n-vector $\Lambda$ on $M$. Note that $\Lambda$ is a Nambu-Poisson structure 
on $M$.\\ 
Recall that a {\bf Nambu-Poisson structure} on $M$ of order $r$ is a skew-symmetric $r$-linear 
map $\{\,,\hdots,\,\}$
$${\cal C}^\infty(M)\times\hdots\times{\cal C}^\infty(M) \longrightarrow {\cal C}^\infty(M),
\quad (f_1,\hdots ,f_r)\longmapsto
\{f_1,\hdots ,f_r\},$$
which satifies
$$\{f_1,\hdots,f_{r-1},gh\}=\{f_1,\hdots,f_{r-1},g\}h+g\{f_1,\hdots,f_{r-1},h\}\quad (L)$$
$$\{f_1,\hdots,f_{r-1},\{g_1,\hdots,g_{r}\}\}=\sum_{i=1}^r
\{g_1,\hdots,g_{i-1},\{f_1,\hdots,f_{r-1},g_i\},g_{i+1},\hdots ,g_r\}\quad (FI)$$
for any $f_1,\hdots ,f_{r-1},g,h,g_1,\hdots ,g_r$ in ${\cal C}^\infty(M)$. 
It is clear that we can associate to such a bracket a $r$-vector on $M$. If $r=2$, we 
rediscover Poisson structures. Thus, Nambu-Poisson structures can be seen as a kind of
generalization of Poisson structures. The notion of Nambu-Poisson structures was introduced
in [T] by Takhtajan in order to give a formalism to an idea of Y. Nambu ([Na]).\\

\noindent Here, we suppose that the set $\{x\in M\,;\, \Lambda_x\neq 0\}$ is dense in $M$.
We are going to associate a cohomology to $(M,\Lambda)$.

\subsection{The choice of the cohomology}
\indent If $M$ is a differentiable manifold of dimension 2, then the Poisson structures on $M$
are the 2-vectors on $M$. If $\Pi$ is a Poisson structure on $M$, then we can associate to
$(M,\Pi)$ the complex
$$0\longrightarrow {\cal C}^\infty(M)\stackrel{\partial}\longrightarrow {\cal X}^1(M)
\stackrel{\partial}\longrightarrow {\cal X}^2(M)\longrightarrow 0$$
with $\partial (g)=[g,\Pi]=X_g$ (Hamiltonian of $g$) if $g\in {\cal C}^\infty(M)$
and $\partial (X)=[X,\Pi]$  ($[\,,\,]$ indicates Schouten's bracket) if $X\in {\cal X}^1(M)$.
The cohomology of this 
complex is called the Poisson cohomology of $(M,\Pi)$. This cohomology has been studied for
instance in [Mo], [N1] and [V].\\

Now if $M$ is of dimension $n$ with $n\geq 3$, we want to generalize this cohomology.
Our first approach was to consider the complex
$$0\longrightarrow \big( {\cal C}^\infty(M)\big) ^{n-1}\stackrel{\partial}\longrightarrow
{\cal X}^1(M)\stackrel{\partial}\longrightarrow {\cal X}^n(M)\longrightarrow 0$$
with $\partial (X)=[X,\Lambda]$ and 
$\partial (g_1,\hdots,g_{n-1})=i_{dg_1\wedge\hdots\wedge dg_{n-1}}\Lambda=
X_{g_1,\hdots,g_{n-1}}$ (Hamiltonian vector field) where we adopt the convention
$i_{dg_1\wedge\hdots\wedge dg_{n-1}}\Lambda=\Lambda(dg_1,\hdots,dg_{n-1},\bullet )$.
We denote $H_\Lambda ^0(M)$, $H_\Lambda ^1(M)$ and $H_\Lambda ^2(M)$ the three spaces of 
cohomology of this complex. 
With this cohomology, we rediscover the interpretation of the first
spaces of the Poisson cohomology, i.e. $H^2_\Lambda(M)$ describes the infinitesimal 
deformations of $\Lambda$ and $H^1_\Lambda(M)$ is the quotient of the algebra of vector 
fields which preserve $\Lambda$ by the ideal of Hamiltonian vector fields.\\

In [I2], the authors associate to any Nambu-Poisson structure on $M$ a cohomology.
The second idea is then to adapt their construction to our particular case.\\
Let $\#_{\Lambda}$ be the morphism of ${\cal C}^\infty(M)$-modules
$\Omega^{n-1}(M)\longrightarrow {\cal X}^1(M)$ :\, $\alpha\mapsto i_{\alpha}\Lambda$.
Note that $\ker \#_\Lambda=\{0\}$ (because the set of regular points of $\Lambda$ is dense).
We can define (see [I1]) a $\R$-bilinear operator $\lcf\; ,\:\rcf :\Omega^{n-1}(M)\times
\Omega^{n-1}(M)\longrightarrow\Omega^{n-1}(M)$ by 
$$\lcf\alpha,\beta\rcf={\cal L}_{\#_\Lambda\alpha}\beta+(-1)^n(i_{d\alpha}
\Lambda)\beta\quad .$$
The vector space $\Omega^{n-1}(M)$ equiped with $\lcf\; ,\: \rcf$ is a Lie algebra (for any 
Nambu-Poisson structure, it is a Leibniz algebra). Moreover
this bracket verifies $\#_\Lambda \lcf\alpha,\beta\rcf=[\#_\Lambda\alpha,\#_\Lambda\beta]$ for 
any $\alpha,\beta$ in $\Omega^{n-1}(M)$. The triple $(\Lambda^{n-1}(T^*(M)),\lcf\; ,\: \rcf  ,
\#_\Lambda)$ is then a Lie algebroid and the Nambu-Poisson cohomology of $(M,\Lambda)$ is the 
Lie algebroid cohomology of $(\Lambda^{n-1}(T^*(M))$ (for any Nambu-Poisson structure, it is more
elaborate see [I2]). More precisely, for every 
$k\in\{0,\hdots,n\}$, we consider the vector space $C^k( \Omega^{n-1}(M);{\cal C}^\infty(M))$
of the skew-symmetric and ${\cal C}^\infty(M)$-$k$-multilinear maps 
$\Omega^{n-1}(M)\times\hdots\times\Omega^{n-1}(M)\longrightarrow{\cal C}^\infty(M)$.
The cohomology operator \\$\partial:C^k( \Omega^{n-1}(M);{\cal C}^\infty(M))\longrightarrow
C^{k+1}( \Omega^{n-1}(M);{\cal C}^\infty(M))$ is defined by
\begin{eqnarray*}
\partial c(\alpha_0,\hdots,\alpha_k)&=&\sum_{i=0}^{k}(-1)^i(\#_\Lambda\alpha_i).c\big(
\alpha_0,\hdots,{\hat {\alpha_i}},\hdots,\alpha_k\big)\\
& &+ \sum_{0\leq i<j\leq k}(-1)^{i+j}
c(\lcf\alpha_i,\alpha_j\rcf,\alpha_0,\hdots,{\hat {\alpha_i}},\hdots,{\hat {\alpha_j}},
\hdots,\alpha_k)
\end{eqnarray*}
for all $c\in C^k(\Omega^{n-1}(M);{\cal C}^\infty(M))$ and $\alpha_0,\hdots,\alpha_k$ in 
$\Omega^{n-1}(M)$.\\
The {\bf {Nambu-Poisson cohomology}} of $(M,\Lambda)$, denoted by $H_{NP}^\bullet(M,\Lambda)$, 
is the cohomology of this complex.

\subsection{An equivalent cohomology}
So defined, the Nambu-Poisson cohomology is quite difficult to manipulate. We are going to give
an equivalent cohomology which is more accessible.\\
Recall that we assume that $M$ admits a volume form $\omega$.\\
Let $f\in{\cal C}^\infty(M)$, we define the operator 
\begin{eqnarray*}
d_f:\Omega^k(M) &\longrightarrow& \Omega^{k+1}(M)\\
         \alpha &\longmapsto&  fd\alpha-kdf\wedge\alpha.
\end{eqnarray*}
It is easy to prove that $d_f\circ d_f=0$. We denote $H_f^\bullet(M)$ the cohomology of this complex.
Let $\flat$ be the isomorphism ${\cal X}^1(M)\longrightarrow\Omega^{n-1}(M)$ \,
$X\longmapsto i_X\omega$.
\begin{lem}
1- If $X\in{\cal X}(M)$, then $\#_\Lambda\big(\flat(X)\big)=(-1)^{n-1}fX$ where 
$f=i_\Lambda\omega$.\\
2- If $X$ and $Y$ are in ${\cal X}(M)$, then $$(-1)^{n-1}\lcf\flat(X),\flat(Y)\rcf=f\flat([X,Y])+
(X.f)\flat(Y)-(Y.f)\flat(X)\,.$$
\end{lem}
{\sf Proof }: 1- Obvious.\\
2-
We have 
$\#_\Lambda\big(\lcf\flat(X),\flat(Y)\rcf\big)=[\#_\Lambda(\flat(X)),\#_\Lambda(\flat(Y))]$
(property of the Lie algebroid), which implies that
\begin{eqnarray*}
\#_\Lambda\big(\lcf\flat(X),\flat(Y)\rcf\big) &=& f(X.f)Y-f(Y.f)X+f^2[X,Y]\\
  &=& (-1)^{n-1}\#_\Lambda\big((X.f)\flat(Y)-(Y.f)\flat(X)+f\flat([X,Y])\big)\,.
\end{eqnarray*}
The result follows via the injectivity of $\#_\Lambda$.$\quad\Box$

\begin{prop}
If we put $f=i_\Lambda\omega$, then $H_{NP}^\bullet(M,\Lambda)$ is isomorphic to 
$H_f^\bullet(M)$.
\end{prop}
{\sf Proof }: For every $k$, we consider the application 
$\varphi: C^k(\Omega^{n-1}(M);{\cal C}^\infty(M))\longrightarrow\Omega^k(M)$ defined by
$$\varphi(c)\big(X_1,\hdots,X_k\big)=c\big((-1)^{n-1}\flat(X_1),\hdots,
(-1)^{n-1}\flat(X_k)\big)\,,$$
where $c\in C^k(\Omega^{n-1}(M);{\cal C}^\infty(M))$ and $X_1,\hdots,X_k\in{\cal X}(M)$. It is 
easy to see that $\varphi$ is an isomorphism of vector spaces. We show that it is an isomorphism
of complexes.\\
Let $c\in C^k(\Omega^{n-1}(M);{\cal C}^\infty(M))$. We put $\alpha=\varphi(c)$. If $X_0,\hdots,
X_k$ are in ${\cal X}(M)$ then
$\varphi(\partial c)(X_0,\hdots,X_k)=(-1)^{(n-1)(k+1)}\partial c (\flat(X_0),\hdots,\flat(X_k))
=A+B$ where\\
$A=(-1)^{(n-1)(k+1)}\sum_{i=0}^k(-1)^i\#_\Lambda\big(\flat(X_i)\big).c\big(\flat(X_0),\hdots,
{\widehat {\flat(X_i)}},\hdots,\flat(X_k)\big)$\\
$B=(-1)^{(n-1)(k+1)}\sum_{0\leq i<j\leq k} (-1)^{i+j}c\big(\lcf\flat(X_i),\flat(X_j)\rcf,
\flat(X_0),
\hdots,{\widehat {\flat(X_i)}},\hdots,{\widehat {\flat(X_j)}},\hdots,\flat(X_k)\big)$.\\
We have $A=f\sum_{i=0}^k (-1)^i X_i.\alpha(X_0,\hdots,{\hat {X_i}},\hdots,X_k)$ and
\begin{eqnarray*}
B \!&=&\! f\sum_{0\leq i<j\leq k} (-1)^{i+j} \alpha([X_i,X_j],X_0,\hdots,{\hat {X_i}},
\hdots,{\hat {X_j}},\hdots,X_k)\\
  & & +\sum_{0\leq i<j\leq k} (-1)^{i+j} (X_i.f)\alpha(X_j,X_0,\hdots,{\hat {X_i}},
\hdots,{\hat {X_j}},\hdots,X_k)\\
  & & -\sum_{0\leq i<j\leq k} (-1)^{i+j} (X_j.f)\alpha(X_i,X_0,\hdots,{\hat {X_i}},
\hdots,{\hat {X_j}},\hdots,X_k)\\
  &=& \!f\sum_{0\leq i<j\leq k} (-1)^{i+j} \alpha([X_i,X_j],X_0,\hdots,{\hat {X_i}},
\hdots,{\hat {X_j}},\hdots,X_k)\\
  & & -k \sum_{i=0}^k (-1)^i (X_i.f)\alpha(X_0,\hdots,{\hat {X_i}},\hdots,X_k)\,.
\end{eqnarray*}
Consequently, $\varphi(\partial c)=d_f\alpha=d_f\big( \varphi(c)\big)$.$\quad\Box$

\begin{rem}
{\rm We claim that this cohomology is a ``good'' generalisation of the Poisson cohomology 
of a 2-dimensional Poisson manifold. Indeed, if $(M,\Pi)$ is an orientable Poisson manifold 
of dimension 2, we consider the volume form $\omega$ on $M$ and we put 
$$\phi^2 : {\cal X}^2(M)\longrightarrow \Omega^2(M) \; {\mbox { and }} \; 
\phi^1 : {\cal X}^1(M)\longrightarrow \Omega^1(M)$$ 
defined by
$$\phi^2(\Gamma)=(i_\Gamma\omega)\omega \quad {\mbox { and }} \quad
\phi^1(X)=-i_X\omega$$
for every 2-vector $\Gamma$ and vector field $X$.\\ 
We also put $\phi^0=id:{\cal C}^\infty(M)\longrightarrow {\cal C}^\infty(M)$.\\
If we denote $\partial$ the operator of the Poisson cohomology, and $f=i_\Pi\omega$, 
it is quite easy to see that
$$\phi:({\cal X}^\bullet(M),\partial)\longrightarrow (\Omega^\bullet(M),d_f)$$
is an isomorphism of complexes.}
\end{rem}

\begin{rems}
{\rm  1- The definitions we have given make sense if we work in the 
holomorphic case or in the formal case.}\\
{\rm 2- {\bf { Important :}}
If $h$ is a function on $M$ which doesn't vanish on $M$, then the cohomologies
$H_f^\bullet(M)$ and $H_{fh}^\bullet(M)$ are isomorphic. \\
Indeed, the applications
$\big(\Omega^k(M),d_{fh}\big) \longrightarrow \big(\Omega^k(M),d_f\big)\quad \alpha\longmapsto
\frac{\alpha}{h^k}$ give an isomorphism of complexes.\\
In particular, if $f$ doesn't vanish on $M$ then $H_f^\bullet(M)$ is isomorphic to the de Rham's
cohomology. }
\label{R1}
\end{rems}

\subsection{Other cohomologies}
We can construct other complexes which look like 
$(\Omega^\bullet(M),d_f)$. More precisely we denote, for $p\in\Z$, 
\begin{eqnarray*}
d_f^{(p)}:\Omega^k(M) &\longrightarrow& \Omega^{k+1}(M)\\
         \alpha &\longmapsto&  fd\alpha-(k-p)df\wedge\alpha.
\end{eqnarray*}
We will denote $H_{f,p}^\bullet(M)$ the cohomology of these complexes. We will see in the next 
section some relations between these different cohomologies.\\

Using the contraction $i_\bullet\omega$, it is quite easy to prove the following proposition. 
\begin{prop}
The spaces $H_\Lambda^1(M)$ and $H_\Lambda^2(M)$ are isomorphic to $H_{f,n-2}^{n-1}(M)$ 
and $H_{f,n-2}^n(M)$.
\end{prop}

\begin{rem}
{\rm The two properties of remark \ref{R1} are valid for $H_{f,p}^\bullet(M)$ with $p\in\Z$.}
\label{R2}
\end{rem} 

\section{Computation}
Henceforth, we will work {\bf {locally}}. Let $\Lambda$ be a germ of $n$-vectors on $\K^n$
($\K$ indicates $\R$ or $\C$ ) with {\underline {$n\geq 3$}}. We denote $\fct$ 
($\Omega^k(\K^n),\chp$) the space of {\bf germs} at 0 of
(holomorphic, analytic, ${\cal C}^\infty$, formal) functions ($k$-forms, vector fields).
We can write $\Lambda$ (with coordinates $(x_1,\hdots,x_n)$) 
$\Lambda=f\frac{\partial }{\partial x_1}\wedge\hdots\wedge \frac{\partial }{\partial x_n}$
where $f\in\fct$. We assume that the volume form $\omega$ is $dx_1\wedge\hdots\wedge dx_n$.\\
We suppose that $f(0)=0$ (see remark \ref{R1}) and that $f$ is of {\bf {finite codimension}}, 
which means that $Q_f=\fct /I_f$ ($I_f$ is the 
ideal spanned by $\frac{\partial f}{\partial x_1},\hdots,\frac{\partial f}{\partial x_n}$) is
a  finite dimensional vector space.
\begin{rem}
{\rm It is important to note that, according to Tougeron's theorem (see for instance
[AGV]), if $f$ is of finite codimension, then the set $f^{-1}(\{0\})$ is, from the topological 
point of view, the same as the set of the zeros of a polynomial.\\
Therefore, if $g$ is a germ at 0 of functions which satisfies $fg=0$, then $g=0$.}
\label{codim}
\end{rem}
Moreover we suppose that $f$ is a {\bf {quasihomogeneous}} polynomial of degree N 
(for a justification of this additional assumption, see section 3). 
We are going to recall the definition of the quasi-homogeneity.

\subsection{Quasi-homogeneity}
Let $(w_1,\hdots,w_n)\in{(\N^*)}^n$. We denote $W$ the vector field 
$w_1x_1\frac{\partial }{\partial x_1}+\hdots +w_nx_n\frac{\partial }{\partial x_n}$ on $\K^n$.
We will say that a tensor $T$ is quasihomogeneous with weights $w_1,\hdots, w_n$ and of 
(quasi)degree $\degp\in\Z$ if ${\cal L}_WT=\degp T$ (${\cal L}$ indicates the Lie derivative
operator). Note that $T$ is then polynomial.\\ 
If $f$ is a \qh polynomial of degree N then ${\rm N}=k_1w_1+\hdots +k_nw_n$ with 
$k_1,\hdots,k_n\in\N$ ; so, an integer is not necessarily the quasidegree of a polynomial.
If $f\in\K\big[[x_1,\hdots,x_n]\big]$, we can write $f=\sum_{i=0}^{\infty}f_i$ with 
$f_i$ \qh  of degree $i$ (we adopt the convention that $f_i=0$ if $i$ is not a
quasidegree); $f$ is said to be of order d (${\mbox {ord}}(f)={\rm d}$) if all of its 
monomials have a degree d or higher. For more details, consult [AGV].\\
Since ${\cal L}_W$ and the exterior differentiation $d$ commute, if $\alpha$ is a \qh
$k$-form, then $d\alpha$ is a \qh $(k+1)$-form of degree $\deg \alpha$. In particular,
it is important to notice that $dx_i$ is a \qh 1-form
of degree $w_i$ (note that $\frac{\partial }{\partial x_i}$ is a \qh vector 
field of degree $-w_i$). Thus, the volume form 
$\omega=dx_1\wedge\hdots\wedge dx_n$ is \qh of degree $w_1+\hdots+w_n$. Note that a \qh 
non zero $k$-form ($k\geq 1$) has a degree strictly positive.\\
Note that if $f$ is a \qh polynomial of degree N, then the $n$-vector 
$\Lambda=f\frac{\partial }{\partial x_1}\wedge\hdots\wedge\frac{\partial }{\partial x_n}$ is \qh
of degree ${\rm N}-\sum_iw_i$.\\
In the sequel, the degrees will be quasidegrees with respect to 
$W=w_1x_1\frac{\partial}{\partial x_1}+\hdots+w_nx_n\frac{\partial}{\partial x_n}$.\\
We will need the following result.
\begin{lem}
Let $k_1,\hdots,k_n\in\N$ and put $p=\sum k_iw_i$. Assume that $g\in\fct$ and 
$\alpha\in\Omega^i(\K^n)$ verify ${\rm {ord}}(j_0^\infty(g))>p$ and 
${\rm {ord}}(j_0^\infty(\alpha))>p$ ($j_0^\infty$ indicates the $\infty$-jet at 0). Then\\
\indent 1. there exists $h\in\fct$ such that $W.h-ph=g$,\\ 
\indent 2. there exists $\beta\in\Omega^i(\K^n)$ such that ${\cal L}_W\beta-p\beta=\alpha$.
\label{Rouss}
\end{lem}
{\sf {proof }}: The first claim is only a generalisation of a lemma given (in dimension 2) 
in [Mo] and it can be proved in the same way. The second claim is a consequence of the first.\\

Now we are going to compute the spaces $H^k_f(\K^n)$ (i.e $H^k_{NP}(\K^n,\Lambda)$) for 
$k=0,\hdots,n$. We will denote
$Z_f^k(\K^n)$ and $B_f^k(\K^n)$ the spaces of $k$-cocycles and $k$-cobords. We will also compute 
some spaces $H_{f,p}^k(\K^n)$ with particular interest in the spaces $H^n_{f,n-2}(\K^n)$ 
(i.e. $H^2_\Lambda(\K^n)$) and $H^{n-1}_{f,n-2}(\K^n)$ (i.e. $H^1_\Lambda(\K^n)$). We will 
denote $Z^k_{f,p}(\K^n)$ \big($B^k_{f,p}(\K^n)$\big) the spaces of $k$-cocycles ($k$-cobords) 
for the operator $d_f^{(p)}$.

\subsection{Two useful preliminary results}
In the computation of these spaces of cohomology, we will need the two following propositions.
The first is only a corollary of the de Rham's division lemma (see [dR]).
\begin{prop}
Let $f\in\fct$ of finite codimension. If $\alpha\in\Omega^k(\K^n)$ $(1\leq k\leq n-1)$
verifies $df\wedge\alpha=0$ then there exists $\beta\in\Omega^{k-1}(\K^n)$ such that
$\alpha=df\wedge\beta$.
\label{div}
\end{prop}

\begin{prop}
Let $f\in\fct$ of finite codimension.
Let $\alpha$ be a k-form $(2\leq k \leq n-1)$ which verifies $d\alpha=0$ and $df\wedge\alpha=0$
then there exists $\gamma\in\Omega^{k-2}(\K^n)$ such that $\alpha=df\wedge d\gamma$.
\label{Vey}
\end{prop}
{\sf {Proof }}: We are going to prove this result in the formal case and in the analytical case.\\
{\it {Formal case}}: Let $\alpha$ be a \qh k-form of degree $p$ which verifies the hypotheses.
Since $df\wedge \alpha=0$, we have $\alpha=df\wedge\beta_1$ where $\beta_1$ is a \qh $(k-1)$-form of
degree $p-\degp$. Now, since $d\alpha=0$, we have $df\wedge d\beta_1=0$ and so 
$d\beta_1=df\wedge\beta_2$, where $\beta_2$ is a \qh $(k-1)$-form of degree $p-2\degp$. This way, we 
can construct a sequence $(\beta_i)$ of \qh $(k-1)$-forms with $\deg \beta_i=p-i\degp$ which 
verifies $d\beta_i=df\wedge\beta_{i+1}$. Let $q\in\N$ such that $p-q\degp\leq 0$. Thus, we have
$\beta_q=0$ and so $d\beta_{q-1}=0$ i.e. $\beta_{q-1}=d\gamma_{q-1}$ where $\gamma_{q-1}$ is a
$(k-2)$-form. Consequently, $d\beta_{q-2}=df\wedge d\gamma_{q-1}$ which implies that 
$\beta_{q-2}=-df\wedge\gamma_{q-1}+d\gamma_{q-2}$, where $\gamma_{q-2}$ is a $(k-2)$-form.
In the same way, $d\beta_{q-3}=df\wedge d\gamma_{q-2}$ so $\beta_{q-3}=-df\wedge \gamma_{q-2}
+d\gamma_{q-3}$ where $\gamma_{q-3}$ is a $(k-2)$-form. This way, we can show that 
$\beta_1=-df\wedge\gamma_2+d\gamma_1$ where $\gamma_1$ and $\gamma_2$ are $(k-2)$-forms.
Therefore, $\alpha=df\wedge d\gamma_1$\\
{\it {Analytical case }}: In [Ma], Malgrange gives a result on the relative cohomology of a germ
of an analytical function. In particular, he shows that in our case, if $\beta$ is a germ at 0
of analytical $r$-forms ($r<n-1$) which verifies $d\beta=df\wedge\mu$ ($\mu$ is a $r$-form) 
then there exists two germs of analytical $(r-1)$-forms $\gamma$ and $\nu$ such that 
$\beta=d\gamma+df\wedge\nu$.\\ 
Now, we are going to prove our proposition. 
Let $\alpha$ be a germ of analytical $k$-forms ($2\leq k\leq n-1$) 
which verifies the hypotheses of the proposition.
Then there exists a $(k-1)$-form $\beta$ such that $\alpha=df\wedge \beta$ 
(proposition \ref{div}). But since 
$0=d\alpha=-df\wedge d\beta$, we have $d\beta=df\wedge\mu$ and so ([Ma])
$\beta=d\gamma+df\wedge\nu$ where $\gamma$ and $\nu$ are analytical $(k-2)$-forms. We deduce that 
$\alpha=df\wedge d\gamma$ where $\gamma$ is analytic. $\quad\Box$

\begin{rem}
{\rm {\bf Important:}\\
In fact, some results which appear in [R] lead us to think that this proposition is not true
in the real ${\cal C}^\infty$ case.\\
The computation of the spaces $H^n_{f,p}(\K^n)$, $H^{n-1}_{f,p}(\K^n)$ ($p\neq n-2$) and 
$H^0_{f,p}(\K^n)$ doesn't use this proposition so, it still holds in the ${\cal C}^\infty$ case.\\
The results we find on the other spaces should be the same in the ${\cal C}^\infty$ case as in
the analytical case but another proof need to be found.}
\end{rem}

\subsection{Computation of $H_{f,p}^0(\K^n)$}
We consider the application $d_f^{(p)}:\Omega^0(\K^n)\longrightarrow\Omega^1(\K^n)\quad
g\longmapsto fdg+pdf\wedge g$.
\begin{theo}
1- If $p>0$ then $H_{f,p}^0(\K^n)=\{0\}$\\
2- If $p\leq 0$ then $H_{f,p}^0(\K^n)=\K {\bf .} f^{-p}$
\label{H0}
\end{theo}
{\sf {Proof }}: 1- If $g\in\fct$ verifies $d_f^{(p)}g=0$ then $d(f^pg)=0$, and so $f^pg$ is
constant. But as $f(0)=0$, $f^pg$ must be 0 i.e. $g=0$ (because $f$ is of finite 
codimension; see remark \ref{codim}).\\
2- We will use an induction to show that for any $k\geq 0$, if $g$ satisfies  $fdg=kgdf$
then $g=\lambda f^k$ where $\lambda\in\K$.\\
For $k=0$ it is obvious.\\
Now we suppose that the property is true for $k\geq 0$. We show that it is still valid for $k+1$.
Let $g\in\fct$ be such that $fdg=(k+1)gdf\, (*)$. Then $df\wedge dg=0$ and so there exists
$h\in\fct$ such that $dg=hdf$ (proposition \ref{div}). Replacing $dg$ by $hdf$ in $(*)$, we get
$fhdf=(k+1)gdf$ i.e. $g=\frac{1}{k+1} fh$. Now, this former relation gives on the one hand
$fdg=\frac{1}{k+1}(f^2dh+fhdf)$ and on the other hand, using $(*)$, $fdg=fhdf$. Consequently,
$fdh=khdf$ and so $h=\lambda f^k$ with $\lambda\in\K$.$\quad\Box$

\subsection{Computation of $H_f^k(\K^n)$ $1\leq k\leq n-2$}
\begin{lem}
Let $\alpha\in Z^k_{f,p}(\K^n)$ with $1\leq k\leq n-2$. Then $\alpha$ is cohomologous to
a closed $k$-form.
\label{LL}
\end{lem}
{\sf {Proof }}: We have $fd\alpha-(k-p)df\wedge\alpha=0$. If $k=p$ then $\alpha$ is closed.
Now we suppose that $k\neq p$.
We put $\beta=d\alpha\in\Omega^{k+1}(\K^n)$. We have
$$0=df\wedge(fd\alpha-(k-p)df\wedge\alpha)=fdf\wedge\alpha$$
so, $df\wedge\alpha=0$.\\
Now, since $d\beta=0$ and $df\wedge\beta=0$, 
proposition \ref{Vey} gives $\beta=df\wedge d\gamma$ with $\gamma\in\Omega^{k-1}(\K^n)$. 
Then, if we consider 
$\alpha^\prime=\alpha-\frac{1}{k-p}\big(fd\gamma-(k-p-1)df\wedge\gamma\big)$, we have 
$d\alpha^\prime=0$ and $fd\gamma-(k-p-1)df\wedge\gamma\in B^k_{f,p}(\K^n)$.$\quad\Box$

\begin{theo} 
If $k\in\{2,\hdots,n-2\}$ then $H_f^k(\K^n)=\{0\}$.
\label{Hk}
\end{theo}
{\sf Proof }: Let $\alpha\in Z_f^k(\K^n)$. Then $\alpha\in\Omega^k(\K^n)$ and verifies 
$fd\alpha-kdf\wedge\alpha=0$.\\
According to lemma \ref{LL} we can assume that $\alpha$ is closed.
Now we show that $\alpha\in B_f^k(\K^n)$.\\
Since $d\alpha=0$ and $df\wedge\alpha=0$, there exists 
$\beta\in\Omega^{k-2}(\K^n)$ such that $\alpha=df\wedge d\beta$ (proposition \ref{Vey}). Thus,
$\alpha=d_f\big(\frac{-1}{k-1}d\beta\big)$.$\quad\Box$

\begin{rem}
{\rm It is possible to adapt this proof to show that $H_{f,p}^k(\K^n)=\{0\}$
if $k\in\{2,\hdots,n-2\}$ and $p\neq k,k-1$.}
\end{rem}

\begin{lem}
Let $\alpha\in Z_f^1(\K^n)$. If ${\rm {ord}}\big( j_0^\infty(\alpha)\big)>\degp$ then 
$\alpha\in B_f^1(\K^n)$.
\label{L1}
\end{lem}
{\sf {Proof }}: According to lemma \ref{LL}, we can assume that $d\alpha=0$.\\ 
Since $df\wedge\alpha=0$ we have $\alpha=gdf$ (proposition \ref{div}) where 
$g$ is in $\fct$  and verifies ${\rm {ord}}\big( j_0^\infty(g)\big) >0$. 
We show that $f$ divides $g$.\\
Let $\bar{g}\in\fct$ be such that $W.\bar{g}=g$ (lemma \ref{Rouss}); note that  
${\rm {ord}}\big( j_0^\infty(\bar{g})\big) >0$.\\ We have 
${\cal L}_W(df\wedge d\bar{g})=\degp df\wedge d\bar{g}+ df\wedge dg$, and since 
$df\wedge dg=-d\alpha=0$, $df\wedge d\bar{g}$ verifies 
$${\cal L}_W(df\wedge d\bar{g})=\degp df\wedge d\bar{g}\,, $$
which means that $df\wedge d\bar{g}$ is either 0 or \qh of degree N. \\ But since 
${\rm {ord}}\big( j_0^\infty(df\wedge d\bar{g})\big) >\degp$, $df\wedge d\bar{g}$ must be 0.\\
Consequently, there exists $\nu\in\fct$ such that 
$\frac{\partial \bar{g}}{\partial x_i}=\nu\frac{\partial f}{\partial x_i}$ for any $i$.
Thus, $W.\bar{g}=\nu W.f$ and so $g=\nu f$.\\
We deduce that $\alpha=f\beta$ with $\beta\in\Omega^1(\K^n)$. \\ Now, we have
$$0=d\alpha=df\wedge\beta +fd\beta\quad {\mbox { and }}\quad 0=df\wedge\alpha=fdf\wedge\beta\,,$$
which implies that $d\beta=0$.\\
Therefore, $\alpha=fdh=d_f(h)$ with $h\in\fct$.$\quad\Box$

\begin{theo}
The space $H_f^1(\K^n)$ is of dimension 1 and spanned by $df$.
\label{H1f}
\end{theo}
{\sf {Proof }}: Let $\alpha\in Z_f^1(\K^n)$. According to lemma \ref{L1} we only have to study 
the case where $\alpha$ is quasihomogeneous with $\deg (\alpha)\leq\degp$. We have
$fd\alpha-df\wedge\alpha=0$ so, $df\wedge d\alpha=0$. We deduce that $d\alpha=df\wedge\beta$
where $\beta$ is a \qh 1-form of degree $\deg (\alpha)-\degp\leq 0$. But since $dx_i$ is
\qh of degree $w_i>0$ for any $i$, every \qh non zero 1-form has a strictly positive 
degree. We deduce that $\beta=0$ and so $d\alpha=0$. Therefore, $df\wedge\alpha=0$ which 
implies that $\alpha=gdf$ where $g$ is a \qh function of degree $\deg (\alpha)-\degp\leq 0$.
Consequently, if $\deg (\alpha)<\degp$ then $g=0$; otherwise $g$ is constant. To conclude, 
note that $df$ is not a cobord because $f$ doesn't divide $df$.\:$\Box$

\subsection{Computation of $H_{f,p}^n(\K^n)$}
We are going to compute the spaces $H_{f,p}^n(\K^n)$ for $p\neq n-1$. We consider the application
$$d_f^{(n-q)}: \Omega^{n-1}(\K^n)\longrightarrow \Omega^{n}(\K^n)\quad 
\alpha \longmapsto fd\alpha-(q-1)df\wedge\alpha$$ with {\underline {$q\neq 1$}} (note that if 
$q=n$ then we 
obtain the space $H^n_{NP}(M,\Lambda)$ and if $q=2$ then we have $H^2_\Lambda(\K^n)$).\\
We denote ${\cal I}^n=\{df\wedge\alpha;\alpha\in\Omega^{n-1}(\K^n)\}$. It is clear that 
${\cal I}^n\simeq I_f$ (recall that $I_f$ is the ideal of $\fct$ spanned by 
$\frac{\partial f}{\partial x_1},\hdots,\frac{\partial f}{\partial x_1}$) and that 
$\Omega^n(\K^n)/{\cal I}^n\simeq Q_f={\cal F}(\K^n)/I_f$.\\
We put $\sigma=i_W\omega$ (recall that $W=w_1x_1\frac{\partial }{\partial x_1}+\hdots+
w_nx_n\frac{\partial }{\partial x_n}$ and that $\omega$ is the standard volume form on $\K^n$).
Note that $\sigma$ is a \qh $(n-1)$-form of degree $\sum_i w_i$ and that 
$dg\wedge\sigma=(W.g)\omega$ if $g\in\fct$.\\
If $\alpha\in\Omega^{n-1}(\K^n)$, we will use the
notation {\underline {${\rm {div}}(\alpha)$}} for $d\alpha={\rm {div}}(\alpha)\omega$;
for example, ${\rm {div}}(\sigma)=\sum_i w_i$. Note that if $\alpha$ is \qh then 
${\rm {div}}(\alpha)$ is \qh of degree $\deg\alpha-\sum_i w_i$.

\begin{lem}
1- If the $\infty$-jet at 0 of $\gamma$ doesn't contain a component of degree $q\degp$ (in
particular if $q\leq 0$) then 
$\gamma\in B^n_{f,n-q}(\K^n)\Leftrightarrow\gamma\in {\cal I}^n$.\\ 
2- If $\gamma$ is a \qh n-form of degree $q\degp$ then 
$\gamma\in B^n_{f,n-q}(\K^n)\Rightarrow\gamma\in {\cal I}^n$.
\label{L2}
\end{lem}
{\sf {Proof }}: If $\gamma=fd\alpha-(q-1)df\wedge\alpha\in B^n\big( d_f^{(n-q)}\big)$
 where $\alpha\in\Omega^{n-1}$ then $\gamma=df\wedge\beta$ with 
$\beta=-(q-1)\alpha+\frac{{\rm {div}}(\alpha)}{\degp}\sigma$.
This shows the second claim and the first part of the first one.\\
Now we prove the reverse of the first claim.\\
{\it {Formal case }}: Let $\gamma=\sum_{i>0} \gamma^{(i)}$ and $\beta=\sum \beta^{(i-\degp)}$
(with $\gamma^{(i)}$ of degree $i$, $\gamma^{(q\degp)}=0$ and $\beta^{(i-\degp)}$ of degree 
$i-\degp$) such that
$\gamma=df\wedge\beta$. If we put $\alpha=\frac{-1}{q-1}\beta
+\sum_i \frac{{\rm {div}}(\beta^{(i-\degp)})}{(q-1)(i-q\degp)}\sigma$, we have
$d_f^{(n-q)}(\alpha)=\gamma$.\\
{\it {Analytical case }}: If $\beta$ is analytic at 0, the function ${\rm {div}}(\beta)$ 
is analytic too
and since $\lim_{i\rightarrow +\infty} \frac{1}{i-q\degp}=0$, the (n-1)-form defined above
is also analytic at 0.\\
{\it {${\cal C}^\infty$ case }}: We suppose that $\gamma=df\wedge\beta$.
If we denote $\tilde{\gamma}=j_0^\infty\big( \gamma\big)$ then
there exists a formal (n-1)-form $\tilde{\alpha}$ such that 
$\tilde{\gamma}=fd\tilde{\alpha}-(q-1)df\wedge\tilde{\alpha}$. Let $\alpha$
be a ${\cal C}^\infty$-(n-1)-form such that $\tilde{\alpha}=j_0^\infty(\alpha)$. This form
verifies $fd\alpha-(q-1)df\wedge\alpha=\gamma+\varepsilon$ where 
$\varepsilon$ is flat at 0. Since $B^n_{f,n-q}(\K^n) \subset {\cal I}^n$, 
$\varepsilon\in {\cal I}^n$ so that $\varepsilon=df\wedge\mu$ where $\mu$ is flat at 0. 
Let $g\in\fct$ be such that $W.g-\big( (q-1)\degp-\sum w_i\big) g=\frac{{\rm {div}}(\mu)}{q-1}$
(lemma \ref{div}).
Then the form $\theta=\frac{-1}{q-1}\mu+g\sigma$ verifies $d_f^{(n-q)}(\theta)=\varepsilon$.
$\quad\Box$

\begin{rem}
{\rm 1- This lemma gives $B^n_{f,n-q}(\K^n)\subset {\cal I}^n$. Thus, there is a surjection from
$H^n_{f,n-q}(\K^n)$ onto $Q_f$. Therefore, if $f$ is not of finite codimension then 
$H^n_{f,n-q}(\K^n)$ is a infinite-dimensional vector space.\\
2- According to this lemma, if $\gamma$ is in ${\cal I}^n$ then there exits a \qh $n$-form 
$\theta$, of degree $q\degp$, such that $\gamma+\theta\in B^n_{f,n-q}(\K^n)$.}
\end{rem}
 
\noindent The first claim of this lemma allows us to state the following theorem.
\begin{theo}
If $q\leq 0$ then $H_{f,n-q}^n(\K^n)\simeq Q_f$.
\label{Hn-}
\end{theo}
Now we suppose that {\underline{$q>1$}}.
\begin{lem}
Let $\alpha\in\Omega^k(\K^n)$ and $p\in\Z$. Then $fd_f^{(p)}(\alpha)=d_f^{(p-1)}(f\alpha)$.
\label{L3}
\end{lem}
{\sf {Proof }}: Obvious.
\begin{lem}
1- Let $q>2$. If $\alpha\in\Omega^n(\K^n)$ is \qh of degree $(q-1)\degp$ and verifies
$f\alpha\in B^n_{f,n-q}(\K^n)$ then $\alpha\in B^n_{f,n-q+1}(\K^n)$.\\
2- If $\alpha$ is \qh of degree $\degp$ with $f\alpha\in B^n_{f,n-2}(\K^n)$ then 
$\alpha=0$.
\label{L4}
\end{lem}
{\sf {Proof }}: 1- We suppose that $\alpha=g\omega$ with $g\in\fct$ \qh of degree
$(q-1)\degp-\sum w_i$. We have $fg\omega=fd\beta-(q-1)df\wedge \beta$ where $\beta$ is a \qh 
(n-1)-form of degree $(q-1)\degp$.\\
If we put $\theta=-(q-1)\beta+\frac{{\rm {div}}(\beta) -g}{\degp}\sigma$ then $df\wedge\theta=0$,
and so $\theta=df\wedge\gamma$ where $\gamma$ is a \qh (n-2)-form of degree $(q-2)\degp$.
Consequently $\beta=\frac{-1}{q-1}df\wedge\gamma+\frac{{\rm {div}}(\beta)-g}{(q-1)\degp}\sigma$.
Now, a computation shows that $fd\beta-(q-1)df\wedge \beta=\frac{1}{q-1}fdf\wedge d\gamma$
i.e. $f\alpha=\frac{1}{q-1}fdf\wedge d\gamma$.\\
Therefore, $\alpha=\frac{1}{q-1}df\wedge d\gamma=\frac{1}{q-1}d_f^{(n-q+1)}\big( \frac{-1}{q-2}
d\gamma\big)$.\\
2- As in 1- (with $q=2$), we have $f\alpha=fg\omega=d_f^{(n-2)}(\beta)$ with $\deg g=\degp$
and $\deg \beta=\degp$. We put
$\theta=-\beta+\frac{{\rm {div}}(\beta)-g}{\degp}\sigma$.\\
If $\theta\neq 0$ then $\theta=df\wedge \gamma$ where $\gamma$ is a \qh (n-2)-form of degree 0
which is not possible. So, $\theta=0$ i.e. $\beta=\frac{{\rm {div}}(\beta)-g}{\degp}\sigma$.\\
We deduce that $fd\beta-df\wedge\beta=0$ i.e. $\alpha=0$.$\quad\Box$\\

\noindent Let ${\cal B}$ be a monomial basis of $Q_f$ (for the existence of such a basis,
 see [AGV]).
We denote $r_j$ ($j=2,\hdots,q-1$) the number of monomials of ${\cal B}$ whose degree is 
$j\degp-\sum w_i$ (this number doesn't depend on the choice of ${\cal B}$). 
We also denote $s$ the dimension
of the space of \qh polynomials of degree $\degp-\sum w_i$ and $c$ the codimension of $f$.

\begin{theo}
Let $\alpha\in\Omega^n(\K^n)$. Then there exist unique  polynomials $h_1,\hdots,h_q$
(possibly zero) such that\\
\indent $\bullet h_1$ is \qh of degree $\degp-\sum w_i$,\\
\indent $\bullet h_j (2\leq j\leq q-1)$ is a linear combination of monomials of ${\cal B}$ of 
degree $j\degp-\sum w_i$,\\
\indent $\bullet h_q$ is a linear combination of monomials of ${\cal B}$ and
$$\alpha=(h_q+fh_{q-1}+\hdots +f^{q-1}h_1)\omega\quad {\rm {mod}}\, B^n_{f,n-q}(\K^n)\,.$$
In particular, the dimension of $H^n_{f,n-q}(\K^n)$ is $c+r_{q-1}+\hdots +r_2+s$.
\label{Hn+}
\end{theo}
{\sf {Proof }}: 
{\it {Existence }}: We suppose that $\alpha=g\omega$ with $g\in\fct$.
There exists $h_q$, a linear combination of the monomials of ${\cal B}$,
such that $g=h_q$ ${\rm {mod}}\, I_f$. So, according to lemma \ref{L2} (see the former remark), 
$g\omega=h_q\omega+df\wedge\beta$ ${\rm {mod}}\, B^n_{f,n-q}(\K^n)$
where $\beta$ is a \qh (n-1)-form of degree $(q-1)\degp$.\\ Consequently, 
$g\omega=h_q\omega+\frac{1}{q-1}fd\beta-\frac{1}{q-1}[fd\beta-(q-1)df\wedge\beta]
\: {\rm {mod }} B^n_{f,n-q}(\K^n)$
so, we can write $$g\omega=h_q\omega+fg_{q-1}\omega\quad {\rm {mod}}\, B^n_{f,n-q}(\K^n)$$ with 
$\deg g_{q-1}=(q-1)\degp-\sum w_i$.\\
In the same way, $$g_{q-1}\omega=h_{q-1}\omega+fg_{q-2}\omega\quad {\rm {mod}}\, 
B^n_{f,n-q+1}(\K^n)$$ where $h_{q-1}$ is a linear combination of the monomials of 
${\cal B}$ of degree $(q-1)\degp-\sum w_i$ and $g_{q-2}$ is \qh of degree 
$(q-2)\degp-\sum w_i \hdots$
$$\hdots \qquad \hdots $$
$\hdots$ and $$g_2\omega=h_2\omega+fh_1\omega\quad {\rm {mod}}\, B^n_{f,n-2}(\K^n)$$ 
where $h_2$ is a linear combination of the monomials of ${\cal B}$ of degree $2\degp-\sum w_i$ 
and $h_1$ is \qh of degree $\degp-\sum w_i$.\\
Using lemma \ref{L3}, we get\\
$$\alpha=g\omega=h_q+h_{q-1}+f^2h_{q-2}+\hdots+f^{q-1}h_1\omega\quad {\rm {mod}}\, 
B^n\big( d_f^{(n-q)}\big)\,.$$
{\it {Unicity }}: Let $g=h_q+fh_{q-1}+\hdots+f^{q-1}h_1$ with $h_1,\hdots,h_q$ as in the 
statement of the theorem. We suppose that $g\omega\in B^n_{f,n-q}(\K^n)$. Then 
$g\omega\in {\cal I}^n$ i.e. $g\in I_f$. But since $fh_{q-1}+\hdots+f^{q-1}h_1\in I_f$
(because $f\in I_f$) we have $h_q\in I_f$ and so $h_q=0$.\\ Now , according to lemma \ref{L4},
$(h_{q-1}+fh_{q-2}+\hdots+f^{q-2}h_1)\omega$ is in $B^n_{f,n-q+1}(\K^n)$ and so, in the
same way, $h_{q-1}=0$.\\
This way, we get $h_q=h_{q-1}=\hdots=h_2=0$ and $fh_1\omega\in B^n_{f,n-2}(\K^n)$.
Lemma \ref{L4} gives $h_1=0$.$\quad\Box$ \\

\noindent This theorem allows us to give the dimension of the spaces $H^n_{NP}(\K^n,\Lambda)$
and $H^2_{\Lambda}(\K^n)$.
\begin{cor}
Let $\alpha\in\Omega^n(\K^n)$. Then there exist unique  polynomials $h_1,\hdots,h_n$
(possibly zero) such that\\
\indent $\bullet h_1$ is \qh of degree $\degp-\sum w_i$,\\
\indent $\bullet h_j (2\leq j\leq n-1)$ is a linear combination of monomials of ${\cal B}$ of 
degree $j\degp-\sum w_i$,\\
\indent $\bullet h_n$ is a linear combination of monomials of ${\cal B}$ and
$$\alpha=(h_n+fh_{n-1}+\hdots +f^{n-1}h_1)\omega\quad {\rm {mod}}\, B^n_f(\K^n)\,.$$
In particular, the dimension of $H^n_{NP}(\K^n,\Lambda)$ is $c+r_{n-1}+\hdots +r_2+s$. 
\label{Hfn}
\end{cor}

\begin{cor}
Let $\alpha\in\Omega^n(\K^n)$. Then there exist unique  polynomials $h_1, h_2$ (possibly zero)
such that\\
\indent $\bullet h_1$ is \qh of degree $\degp-\sum w_i$,\\
\indent $\bullet h_2$ is a linear combination of monomials of ${\cal B}$ and
$$\alpha=(h_2+fh_1)\omega\quad {\rm {mod}}\, B^n_{f,n-2}(\K^n)\,.$$
In particular, the dimension of $H^2_{\Lambda}(\K^n)$ is $c+s$. 
\label{HL2}
\end{cor} 

\begin{rem}
{\rm If $q=1$ then the space $H^n_{f,n-1}(\K^n)$ is $\Omega^n(\K^n)/f\Omega^n(\K^n)$ which is of 
infinite dimension.}
\end{rem}

\subsection{Computation of $H_{f,p}^{n-1}(\K^n)$}
We are going to compute the spaces $H_{f,p}^{n-1}(\K^n)$ with $p\neq n-1$. We consider the piece
of complex $$\Omega^{n-2}(\K^n)\longrightarrow\Omega^{n-1}(\K^n)\longrightarrow
\Omega^n(\K^n)$$ \\
with $d_f^{(n-q)}(\alpha)=fd\alpha-(q-2)df\wedge\alpha$ if $\alpha\in\Omega^{n-2}(\K^n)$,\\
and $d_f^{(n-q)}(\alpha)=fd\alpha-(q-1)df\wedge\alpha$ if $\alpha\in\Omega^{n-1}(\K^n)$
with $q\neq 1$.//
Remember that if $q=n$ we obtain $H^{n-1}_{NP}(K^n,\Lambda)$ and if $q=2$ we have
$H^1_\Lambda(\K^n)$.
\begin{lem}
If $\alpha\in Z^{n-1}_{f,n-q}(\K^n)$ then $\alpha=\frac{{\rm {div}}(\alpha)}{(q-1)\degp}
\sigma+df\wedge\beta$ with $\beta\in\Omega^{n-2}(\K^n)$ and so, $d\alpha$ verifies 
${\cal L}_W(d\alpha)-(q-1)\degp d\alpha=(q-1)\degp df\wedge d\beta$.
\label{L5}
\end{lem}
{\sf {Proof }}: It is sufficient to notice that $df\wedge\big( \alpha- 
\frac{{\rm {div}}(\alpha)}{(q-1)\degp}\sigma\big)=0$ (proposition \ref{div}). For the second 
claim, we have 
$(q-1)\degp d\alpha=\big( W.{\rm {div}}(\alpha)+(\sum w_i){\rm {div}}(\alpha)\big)\omega
-(q-1)\degp df\wedge d\beta$ and the conclusion follows.$\quad\Box$

\begin{lem}
If $\alpha\in Z^{n-1}_{f,n-q}(\K^n)$ with ${\rm {ord}}\big( j_0^\infty(\alpha)\big)>
(q-1)\degp$ then $\alpha$ is cohomologous to a closed (n-1)-form. In particular, if $q\leq 0$
then every $(n-1)$-cocycle for $d^{(n-q)}_f$ is cohomologous to a closed $(n-1)$-form.
\label{L6}
\end{lem}
{\sf {Proof }}: We have $\alpha=\frac{{\rm {div}}(\alpha)}{(q-1)\degp}\sigma+df\wedge\beta$
(lemma \ref{L5})  with $${\cal L}_W(d\alpha)-(q-1)\degp d\alpha=
(q-1)\degp df\wedge d\beta\quad (*)\: .$$
Now, let $\gamma\in\Omega^{n-2}(\K^n)$ such that 
${\cal L}_W\gamma-(q-2)\degp\gamma=(q-1)\degp\beta$ ($\gamma$ exists because ${\rm {ord}}
\big( j_0^\infty(\beta)\big) >(q-2)\degp$, see lemma \ref{Rouss}).\\ We have
${\cal L}_W d\gamma-(q-2)\degp d\gamma=(q-1)\degp d\beta$. Thus $df\wedge d\gamma$ verifies
$${\cal L}_W(df\wedge d\gamma)-(q-1)\degp df\wedge d\gamma=(q-1)\degp df\wedge d\beta
\quad (**)\: .$$
From $(*)$ and $(**)$ we get $d\alpha=df\wedge d\gamma$.\\ Indeed, ${\cal L}_W(d\alpha-df\wedge 
d\gamma)=(q-1)\degp (d\alpha-df\wedge d\gamma)$ but $d\alpha-df\wedge d\gamma$ is not \qh of
degree $(q-1)\degp$.\\
Now, if we put $\theta=\alpha-\frac{1}{q-1}\big( fd\gamma-(q-2)df\wedge\gamma\big)$, we have
$d\theta=0$ and $\theta=\alpha\quad {\rm {mod}}\, B^{n-1}_{f,n-q}(\K^n)$.$\Box$\\

\noindent This lemma allows us to state the following theorem.
\begin{theo}
If we suppose that $q\leq 0$ then $H^{n-1}_{f,n-q}(\K^n)=\{0\}$.
\label{Hn-1qneg}
\end{theo}
{\sf Proof }:
Let $\alpha\in Z^{n-1}_{f,n-q}(\K^n)$. We can suppose (according to the former lemma) that 
$d\alpha=0$. Thus we have $df\wedge\alpha=0$. Proposition \ref{Vey} gives then, 
$\alpha=df\wedge d\gamma$ with $\gamma\in\Omega^{n-3}(\K^n)$. 
Therefore, $\alpha=d_f^{(n-q)}\big( -\frac{1}{q-2}d\gamma\big)$.$\quad\Box$\\

\noindent Now, we assume that {\underline {$q>1$}}.
\begin{lem}
If $\alpha\in Z^{n-1}_{f,n-q}(\K^n)$ is a \qh (n-1)-form whose degree is strictly lower than
$(q-1)\degp$ then $\alpha$ is cohomologous to a closed (n-1)-form.
\label{L7}
\end{lem}
{\sf {Proof }}: According to lemma \ref{L5}, we have 
$\alpha=\frac{{\rm {div}}(\alpha)}{(q-1)\degp}\sigma+df\wedge\beta$
and so, $$d\alpha=\frac{(q-1)\degp}{\deg (\alpha)-(q-1)\degp}df\wedge d\beta\,.$$
We deduce that, if we put 
$\theta=\alpha-d_f^{(n-q)}\big( \frac{\degp}{\deg (\alpha)-(q-1)\degp} d\beta\big)$, we have
$d\theta=0$.$\quad\Box$

\begin{rem}
{\rm A consequence of lemmas \ref{L6} and \ref{L7} is that, if $q>1$, every cocycle 
$\alpha\in Z^{n-1}_{f,n-q}(\K^n)$ is cohomologous to a cocycle $\eta+\theta$ where $\eta$ is in
$Z^{n-1}_{f,n-q}(\K^n)$ and is closed, and $\theta$ is \qh of degree  $(q-1)\degp$.}
\label{Rn-1}
\end{rem}

\begin{lem}
Let $\alpha=g\sigma$ where $g$ is a \qh polynomial of degree $(q-1)\degp-\sum w_i$. Then\\
\indent 1- If $q>2$ then, $\alpha\in B^{n-1}_{f,n-q}(\K^n)\Longleftrightarrow g\omega\in
B^n_{f,n-q+1}(\K^n)$.\\
\indent 2- If $q=2$, $\alpha\in B^{n-1}_{f,n-2}(\K^n)\Longleftrightarrow \alpha=0$.
\label{L8}
\end{lem}
{\sf {Proof }}:1- $\bullet$ We suppose that $\alpha\in B^{n-1}_{f,n-q}(\K^n)$
i.e. $\alpha=fd\beta-(q-2)df\wedge\beta$ with
$\beta\in\Omega^{n-2}(\K^n)$. Then $d\alpha=(q-1)df\wedge d\beta$.\\
On the other hand, $d\alpha=(q-1)\degp g\omega$ so 
$g\omega=\frac{1}{\degp}df\wedge d\beta=d_f^{(n-q+1)}\big(-\frac{d\beta}{(q-2)\degp}\big)$.\\
$\bullet$ Now we suppose that $g\omega\in B^n_{f,n-q+1}(\K^n)$ i.e.
$g\omega=fd\beta-(q-2)df\wedge\beta$ where $\beta$ is a \qh (n-1)-form of degree $(q-2)\degp$. 
We put $\gamma=i_W\beta\in\Omega^{n-2}(\K^n)$. We have
\begin{eqnarray*}
d_f^{(n-q)}(\gamma) &=& fd\gamma-(q-2)df\wedge\gamma\\
                    &=& fd(i_W\beta)-(q-2)df\wedge (i_W\beta)\\
                    &=& f\big( {\cal L}_W\beta-i_W d\beta\big)
                       -(q-2)\big[-i_W(df\wedge\beta)+(i_W df)\wedge\beta\big]\\
                    &=& f(q-2)\degp\beta -i_W\big[ fd\beta-(q-2)df\wedge\beta\big]
                        -(q-2)(W.f)\beta\\
                    &=& -i_W\big[ fd\beta-(q-2)df\wedge\beta\big]\,.
\end{eqnarray*}
Consequently, $d_f^{(n-q)}(\gamma)=-i_W(g\omega)=-g\sigma$.\\
2- If $\alpha=fd\beta$ where $\beta$ is a \qh (n-2)-form of degree $\deg \alpha-\degp=0$ 
then $\beta=0$ and so $\alpha=0$.$\quad\Box$\\
 
\noindent We recall that ${\cal B}$ indicates a monomial basis of $Q_f$.
We adopt the same notations as for theorem \ref{Hn+}.

\begin{theo}
We suppose that $q>2$.
Let $\alpha\in Z^{n-1}_{f,n-q}(\K^n)$. There exist unique polynomials 
$h_1,\hdots,h_{q-1}$ (possibly zero) such that\\
\indent $\bullet$ $h_1$ is \qh of degree $\degp-\sum w_i$,\\
\indent $\bullet$ $h_k$ ($k\geq 2$) is a linear combination of monomials of ${\cal B}$ of
degree $k\degp-\sum w_i$ and
$$\omega=(h_{q-1}+fh_{q-2}+\hdots+f^{q-2}h_1)\sigma\quad {\rm {mod}}\, 
B^{n-1}_{f,n-q}(\K^n)\,.$$
In particular, the dimension of the space $H^{n-1}_{f,n-q}(\K^n)$ is 
$r_{q-1}+\hdots +r_2+s$.
\label{Hn-1q}
\end{theo}
{\sf {Proof }}: If $\alpha\in Z^{n-1}_{f,n-q}(\K^n)$ then  $\alpha$ is cohomologous to 
$\eta+\theta$, where $\eta$ is in $Z^{n-1}_{f,n-q}(\K^n)$ and is closed, and $\theta$ is \qh of 
degree  $(q-1)\degp$ (see remark \ref{Rn-1}).\\
The same proof as in theorem \ref{Hn-1qneg} shows that $\eta$ is a cobord.\\
Now, we have to study $\theta$. According to lemma \ref{L5}, we can write
$\theta=\frac{{\rm {div}}(\theta)}{(q-1)\degp}\sigma+df\wedge\beta$  
($\beta\in\Omega^{n-2}(\K^n)$) with 
${\cal L}_W(d\theta)-(q-1)\degp d\theta=(q-1)\degp df\wedge d\beta$. Since $\theta$ is \qh of
degree $(q-1)\degp$, the former relation gives $df\wedge d\beta=0$. Consequently, if we put
$\gamma=df\wedge\beta$, proposition \ref{Vey} gives $\gamma=df\wedge d\xi$.\\
Therefore, $\gamma=d_f^{(n-q)}\big(-\frac{1}{q-2}d\xi\big)$ and so 
$\theta=\frac{{\rm {div}}(\theta)}{(q-1)\degp}\sigma\: {\rm {mod}}\, 
B^{n-1}_{f,n-q}(\K^n)$. The conclusion follows using lemma \ref{L8} and theorem
\ref{Hn+}.$\quad\Box$

\begin{cor}
We suppose that $q=n$. 
Let $\alpha\in Z^{n-1}_f(\K^n)$. There exist unique polynomials 
$h_1,\hdots,h_{n-1}$ (possibly zero) such that\\
\indent $\bullet$ $h_1$ is \qh of degree $\degp-\sum w_i$,\\
\indent $\bullet$ $h_k$ ($k\geq 2$) is a linear combination of monomials of ${\cal B}$ of
degree $k\degp-\sum w_i$ and
$$\omega=(h_{n-1}+fh_{n-2}+\hdots+f^{n-2}h_1)\sigma\quad {\rm {mod}}\, 
B^{n-1}_f(\K^n)\,.$$
In particular, the dimension of the space $H^{n-1}_{NP}(\K^n,\Lambda)$ is 
$r_{n-1}+\hdots +r_2+s$.
\label{Hfn-1}
\end{cor}

\begin{rem}
{\rm If $q=2$, the description of the space $H^{n-1}_{f,n-2}(\K^n)$ (and so $H^1_\Lambda(\K^n)$) 
is more difficult. It is possible
to show that this space is not of finite dimension. Indeed, let us consider the case 
$n=3$ in order to simplify (but it is valid for any $n\geq 3$). We put 
$\alpha=g\big( \frac{\partial f}{\partial x}dx\wedge dz+
\frac{\partial f}{\partial y}dy\wedge dz\big)$ where $g$ is a function which depends only on
$z$. We have $d\alpha=0$ and $df\wedge\alpha=0$ so $\alpha\in Z^{n-1}_{f,n-2}(\K^n)$
but $\alpha\not\in B^n_{f,n-2}(\K^n)$ because $f$ doesn't divide $\alpha$.}
\end{rem} 
We can yet give more precisions on the space $H^{n-1}_{f,n-2}(\K^n)$.

\begin{theo}
Let $E$ be the space of (n-1)-forms $h\sigma$ where $h$ is a \qh polynomial of degree 
$\degp-\sum w_i$, and $F$
the quotient of the vector space $\{ df\wedge d\gamma;\gamma\in\Omega^{n-3}(\K^n)\}$
by the subspace $\{df\wedge d(f\beta);\beta\in\Omega^{n-3}(\K^n)\}$.\\
Then $H^{n-1}_{f,n-2}(\K^n)=E \oplus F$.
\label{HL1}
\end{theo}
{\sf {Proof }}: Let $\alpha$ in $Z^{n-1}_{f,n-2}(\K^n)$.\\
According to remark \ref{Rn-1}, there exixt a closed $(n-1)$-form $\eta$ with 
$\eta\in Z^{n-1}_{f,n-2}(\K^n)$ and a \qh $(n-1)$-form $\theta$, such that $\alpha$ is cohomologous
to $\eta+\theta$.\\
We have (lemma \ref{L5})
$\theta=\frac{{\rm {div}}(\theta)}{\degp}\sigma+df\wedge\beta$ with $\beta$ \qh of degree
0 which is possible only if $\beta\neq 0$. So, $\theta=g\sigma$ where $g$ is a \qh polynomial
of degree $\degp-\sum w_i$. Lemma \ref{L8} says that $\theta\in B^{n-1}_{f,n-2}(\K^n)$
if and only if $\theta=0$.\\
Now we study $\eta$. Proposition \ref{Vey} gives $\eta=df\wedge d\gamma$ where $\gamma$ is
a $(n-3)$-form. If we suppose that $\eta\in B^{n-1}_{f,n-2}(\K^n)$ then 
$df\wedge d\gamma=fd\xi$ with $\xi\in\Omega^{n-2}(K^n)$ and so, $df\wedge d\xi=0$. Now we 
apply proposition \ref{Vey} to $d\xi$ and we obtain $d\xi=df\wedge d\beta$ with 
$\beta\in\Omega^{n-3}(\K^n)$.  
Consequently, $df\wedge d\gamma=fdf\wedge d\beta$ which implies that 
$d\gamma=fd\beta+df\wedge\mu$ with $\mu\in\Omega^{n-3}(\K^n)$, and so
$d\gamma=d(f\beta)+df\wedge\nu$ with $\nu\in\Omega^{n-3}(\K^n)$.\\
Therefore, $\eta\in B^{n-1}_{f,n-2}(\K^n)\Leftrightarrow\eta=df\wedge d(f\beta)$.$\quad\Box$

\subsection{Summary}
It is time to sum up the results we have found.\\
\indent The cohomology $H_f^\bullet(\K^n)$ (and so the Nambu-Poisson cohomology 
$H^\bullet_{NP}(\K^n,\Lambda)$) has been entirely computed
(see theorems \ref{H0}, \ref{Hk}, \ref{H1f}, and corollaries \ref{Hfn} and \ref{Hfn-1}) :\\
The spaces of this cohomology are of finite dimension and only the "extremal" ones (i.e
$H^0$, $H^1$, $H^{n-1}$ and $H^n$) are possibly different to $\{0\}$. The spaces 
$H^0_{NP}(\K^n,\Lambda)$ and $H^1_{NP}(\K^n,\Lambda)$ are always of dimension 1. The dimensions
of the spaces $H^{n-1}_{NP}(\K^n,\Lambda)$ and $H^n_{NP}(\K^n,\Lambda)$ depend on the one hand
on the type of the singularity of $\Lambda$ (via the role played by $Q_f$), and on the other 
hand, on the "polynomial nature" of $\Lambda$.\\
 
Concerning the cohomology $H^\bullet_{f,n-2}(\K^n)$, we have 
computed $H^n$, i.e. $H^n_\Lambda(\K^n)$ (see corollary \ref{HL2}) and we have given a sketch 
of description of $H^{n-1}$ (see theorem \ref{HL1}). We have also computed the spaces 
$H^0_{f,n-2}(\K^n)$ (theorem \ref{H0}) and $H^k_{f,n-2}(\K^n)$ (theorem \ref{Hk}) for 
$k\neq n-2,n-1$, but these spaces are not particularly interesting for our problem.\\
The space $H^2_\Lambda(\K^n)$, which describes the 
infinitesimal deformations of $\Lambda$ is of finite dimension and its dimension has the same 
property as the dimension of $H^n_{NP}(\K^n,\Lambda)$. On the other hand, the space 
$H^1_\Lambda(\K^n)$ which is the space of the vector fields preserving $\Lambda$ modulo
the Hamiltonian vector fields, is not of finite dimension.\\

\noindent It is interesting to compare the results we have found on these two cohomologies 
with the ones given in [Mo] on the computation of the Poisson cohomology in dimension 2.\\
  
Finally, if $p\neq 0, n-2, n-1$ we have computed the spaces $H^0_{f,p}(\K^n)$, 
$H^{n-1}_{f,p}(\K^n)$, $H^n_{f,p}(\K^n)$ and $H^k_{f,p}(\K^n)$ with
$k\neq p, p+1$.\\
If $p=n-1$ we have computed the spaces $H^0_{f,n-1}(\K^n)$ and $H^k_{f,n-1}(\K^n)$ 
for $2\leq k\leq n-2\quad k\neq p, p+1$ (the space $H^n_{f,n-1}(\K^n)$ is of infinite dimension).

\section{Examples}
In this section, we will explicit the cohomology of some particular germs of $n$-vectors.

\subsection{Normal forms of $n$-vectors}
Let $\Lambda=f\frac{\partial }{\partial x_1}\wedge\hdots\wedge\frac{\partial }{\partial x_n}$
be a germ at 0 of n-vectors on $\K^n$ ({\underline {$n\geq 3$}}) with $f$ of finite codimension 
(see the beginning of section 3) and $f(0)=0$ (if $f(0)\neq 0$, then the local triviality theorem, 
see [AlGu], [G] or [N2], allows us to write, up to a change of coordinates, that
$\Lambda=\frac{\partial }{\partial x_1}\wedge\hdots\wedge
\frac{\partial }{\partial x_n}$).
\begin{prop}
If 0 is not a critical point for $f$ then there exist local coordinates $y_1,\hdots,y_n$
such that 
$$\Lambda=y_1\frac{\partial }{\partial y_1}\wedge\hdots\wedge
\frac{\partial }{\partial y_n}\, .$$
\label{reg}
\end{prop}
{\it {Proof }}: A similar proposition is shown for instance in [Mo] in dimension 2. The proof
can be generalized to the $n$-dimensional ($n\geq 3$) case.\\

Now we suppose that 0 is a critical point of $f$. Moreover, we suppose that the germ $f$ is 
{\bf {simple}}, which means that a sufficiently small neighbourhood (with respect to Whitney's 
topology; see [AGV]) of $f$ intersects only a finite number of R-orbits (two germs $g$ and 
$h$ are said R-equivalent if there exits $\varphi$, a local diffeomorphism at 0, such that 
$g=h\circ\varphi$). Simple germs are those who present a certain kind of stability under
deformation.\\
The following theorem can be found in [A] with only sketches of the proofs. In [Mo], a
similar theorem (in dimension 2) is proved and the demonstration can be adapted here.
\begin{theo}
Let $f$ be a simple germ at 0 of finite codimension. Suppose that $f$ has at 0 a critical point
with critical value 0. Then there exist local coordinates $y_1,\hdots,y_n$ such that the germ
$\Lambda=f\frac{\partial }{\partial x_1}\wedge\hdots\wedge\frac{\partial }{\partial x_n}$
can be written, up to a multiplicative constant, 
$g\frac{\partial }{\partial y_1}\wedge\hdots\wedge\frac{\partial }{\partial y_n}$ where $g$ is
in the following list.
\begin{eqnarray*}
A_k &:& y_1^{k+1}\pm y_2^2\pm\hdots\pm y_n^2\quad k\geq 1\\
D_k &:& y_1^2y_2\pm y_2^{k-1}\pm y_3^2\pm\hdots\pm y_n^2\quad k\geq 4\\
E_6 &:& y_1^3+y_2^4\pm y_3^2\pm\hdots\pm y_n^2\\
E_7 &:& y_1^3+y_1y_2^3\pm y_3^2\pm\hdots\pm y_n^2\\
E_8 &:& y_1^3+y_2^5\pm y_3^2\pm\hdots\pm y_n^2
\end{eqnarray*}
\label{fn}
\end{theo}

Proposition \ref{reg} and theorem \ref{fn} describe most of the germs at 0 of $n$-vectors on 
$\K^n$ vanishing at 0.\\ 
We can notice that the models given in the former list are all \qh polynomials; which 
justifies the assumption we made in section 2. 

\subsection{Some examples}
1- The regular case : $f(x_1,\hdots,x_n)=x_1$ .\\
It is easy to see that $Q_f=\{0\}$ and that $f$ is \qh of degree
$\degp=1$, with respect to $w_1=\hdots =w_n=1$. We have $\degp-\sum w_i<0$, so
$H_f^0(\K^n)\simeq\K$ , $H_f^1(\K^n)=\K {\bf .} dx_1$ and $H_f^k(\K^n)=\{0\}$ 
for any $k\geq 2$.\\

\noindent 2- Non degenerate singularity: $f(x_1,\hdots,x_n)=x_1^2+\hdots +x_n^2$ with 
$n\geq 3$.\\
We have $\degp=2$ and $w_1=\hdots =w_n=1$. The space $Q_f$ is isomorphic to $\K$ and
is spanned by the constant germ 1, which is of degree 0. \\
We deduce that $H_f^0(\K^n)\simeq \K$, $H_f^1(\K^n)=\K.(x_1dx_1+\hdots+x_ndx_n)$ and 
$H_f^k=\{0\}$ for $2\leq k\leq n-2$.\\
In order to describe the spaces $H^{n-1}_f(\K^n)$ and $H^n_f(\K^n)$, we look for an integer 
$k\in\{1,\hdots,n-1\}$ such that $k\degp-\sum w_i=\deg 1$ i.e. $2k-n=0$. \\
Therefore,\\
\indent if $n$ is even then $\{\omega,f^{\frac{n}{2}}\omega\}$ is a basis of $H_f^n(\K^n)$ and
$H_f^{n-1}(\K^n)$ is spanned by $\{f^{\frac{n}{2}-1}\sigma\}$\\
\indent if $n$ is odd then $H_f^{n-1}(\K^n)=\{0\}$ and the space $H_f^n(\K^n)$ is spanned by 
$\{\omega\}$.\\
We recall that $\omega=dx_1\wedge\hdots\wedge dx_n$ and 
$$\sigma=i_W\omega=\sum_{i=1}^n (-1)^{i-1}x_idx_1\wedge\hdots\wedge {\widehat {dx_i}}
\wedge\hdots\wedge dx_n\,.$$

\noindent 3- The case $A_2$ with $n=3$: $f(x_1,x_2,x_3)=x_1^3+x_2^2+x_3^2$.\\
Here, $w_1=2$, $w_2=w_3=3$ and $\degp=6$. Thus, $\degp-\sum w_i=-2$, 
$2\degp-\sum w_i=4$ and $3\degp-\sum w_i=10$.\\
Moreover, ${\cal B}=\{1,x_1\}$ is a monomial basis of $Q_f$. But as $\deg 1=0$
and $\deg x_1=3$, we have: $$H_f^0(\K^3)\simeq\K\,, 
H_f^1(\K^3)=\K . (3x_1dx_1+2x_2dx_2+2x_3dx_3)$$ 
$${\mbox {and}} \quad H_f^2(\K^3)=H_f^3(\K^3)=\{0\}\,.$$

\noindent 4- The case $D_5$ with $n=4$: $f(x_1,x_2,x_3,x_4)=x_1^2x_2+x_2^4+x_3^2+x_4^2$.\\
We have $w_1=3$, $w_2=2$, $w_3=w_4=4$ and $\degp=8$ then $\degp-\sum w_i=-5$,
$2\degp-\sum w_i=3$, $3\degp-\sum w_i=11$ and $4\degp-\sum w_i=19$.\\
Now, ${\cal B}=\{1,x_1,x_2,x_2^2,x_2^3\}$ is a monomial basis of $Q_f$. Here, $\deg 1=0$,
$\deg x_1=3$, $\deg x_2=2$, $\deg x_2^2=4$ and $\deg x_2^3=6$. Thus, the only element of 
${\cal B}$ whose degree is of type $k\degp-\sum w_i$ is $x_1$.\\
Consequently, $$H_f^0(\K^4)\simeq\K\,,
H_f^1(\K^4)=\K\ . \big( 2x_1x_2dx_1+(x_1^2+4x_2^3)dx_2+2x_3dx_3+2x_4dx_4\big)\,,$$ 
$$H_f^2(\K^4)=\{0\}\,, H_f^3(\K^4)=\K {\bf .} (x_1\sigma)$$ and 
$\{\omega,x_1\omega,x_2\omega,x_2^2\omega,x_2^3\omega,x_1f\omega\}$ is a basis of 
$H_f^4(\K^4)$.\\
Here, we have $W=3x_1\frac{\partial }{\partial x_1}+
2x_2\frac{\partial }{\partial x_2}+4x_3\frac{\partial }{\partial x_3}+
4x_4\frac{\partial }{\partial x_4}$ and \\ 
$\sigma=3x_1dx_2\wedge dx_3\wedge dx_4 -2x_2dx_1\wedge dx_3\wedge dx_4
+4x_3dx_1\wedge dx_2\wedge dx_4- 4x_4dx_1\wedge dx_2\wedge dx_3$.

\end{document}